\documentclass[article]{amsart}
\usepackage{amssymb}
\usepackage{hyperref}
\numberwithin{equation}{section}
\newtheorem{theorem}{Theorem}[section]
\newtheorem{corollary}[theorem]{Corollary}
\newtheorem{lemma}[theorem]{Lemma}

\begin{document}
\title[Zeros distribution and uniqueness of derivatives of difference]{Some results on zeros distributions and uniqueness of derivatives of difference polynomials}
\author[K Liu, X. L Liu and and T.B Cao]
{Kai Liu, Xinling Liu, Tingbin Cao}

\address{Department of Mathematics\\
Nanchang University\\
Nanchang, Jiangxi, 330031, P.R. China}
 \email{liukai418@126.com}
 \thanks{This work was partially supported by the
NNSF (No. 11026110), the NSF of Jiangxi (No. 2010GQS0144,
2010GQS0139) and the YFED of Jiangxi (No. GJJ11043, No. GJJ10050) of
China.}

\address{Department of Mathematics\\
Nanchang University\\
Nanchang, Jiangxi, 330031, P.R. China}
\email{sdliuxinling@hotmail.com}

 \address{Ting-Bin Cao \newline Department of Mathematics\\
Nanchang University\\
Nanchang, Jiangxi, 330031, P.R. China}
 \email{tbcao@ncu.edu.cn}


\begin{abstract}
We consider the zeros distributions on the derivatives of difference
polynomials of meromorphic functions, and present some results which
can be seen as the discrete analogues of Hayman conjecture
\cite{hayman1}, also partly answer the question given in
\cite[P448]{luolin}. We also investigate the uniqueness problems of
difference-differential polynomials of entire functions sharing one
common value. These theorems improve the results of Luo and
Lin\cite{luolin} and some results of present authors
\cite{liuliucao}.
\end{abstract}

\maketitle

\section{Introduction}
In this paper, a meromorphic function $f$ means meromorphic in the
complex plane. If no poles occur, then $f$ reduces to an entire
function. Throughout of this paper, we denote by $\rho(f)$ and
$\rho_{2}(f)$ the order of $f$ and the hyper order of $f$ \cite{I.L,
yang and yi}. In addition, if $f-a$ and $g-a$ have the same zeros,
then we say that $f$ and $g$ share the value $a$ IM (ignoring
multiplicities). If $f-a$ and $g-a$ have the same zeros with the
same multiplicities, then $f$ and $g$ share the value $a$ CM
(counting multiplicities). We assume that the reader is familiar
with standard symbols and fundamental results of Nevanlinna Theory
\cite{hayman3, I.L, yang and yi}.

\medskip

 Given a meromorphic function $f(z)$, recall that
$\alpha(z)\not\equiv0, \infty$ is a small function with respect to
$f(z)$, if $T(r,\alpha)=S(r,f)$, where $S(r,f)$ is used to denote
any quantity satisfying $S(r,f)=o(T(r,f))$, and $r\rightarrow\infty$
outside of a possible exceptional set of finite logarithmic measure.

\medskip

The following result is related to Hayman conjecture \cite[Theorem
10]{hayman1} which has been considered in several papers later, such
as \cite{ww, chen, mues}.

\medskip

{\noindent\bf Theorem A}. \cite[Theorem 1]{chen} Let $f$ be a
transcendental meromorphic function. If $n\geq1$ is a positive
integer, then $f^{n}f'-1$ has infinitely many zeros.

\medskip

Remark that
 $[f^{n+1}]'=(n+1)f^{n}f'$ in Theorem A, Chen \cite{chen22}, Wang and Fang \cite{wangyuefei,wangyuefei and
Fang} improved Theorem A by proving the following result.

\medskip

{\noindent\bf Theorem B}. Let $f$ be a transcendental entire
function, $n, k$ be two positive integers with $n\geq k+1$. Then
$(f^{n})^{(k)}-1$ has infinitely many zeros.

\medskip

Laine and Yang \cite{I.L and yang} firstly investigated the zeros of
$f(z)^{n}f(z+c)$ and proved the following result.

\medskip

{\noindent\bf Theorem C.} Let $f$ be a transcendental entire
function of finite order and c be nonzero complex constant. If
$n\geq2$, then $f(z)^nf(z+c)-a$ has infinitely many zeros, where
$a\in\mathbb{C}\backslash\{0\}$.

\medskip

Recently, some papers are devoting to improve Theorem C, the
constant $a$ can be replaced by a nonzero polynomial \cite{liukai
arch} or by a small function $a(z)$ \cite{liuliucao}. In addition,
\cite{rocky, liucao,luolin,zhang11} are devoting to the cases of
meromorphic function $f$ or more general difference products. In the
following, without special stated, we assume that $c$ is a nonzero
constant, $n,m,k,s,t$ are positive integers, $a(z)$ is a nonzero
small function with respect to $f(z)$. Let
$P(z)=a_nz^n+a_{n-1}z^{n-1}+\cdots+a_1z+a_0$ be a nonzero
polynomial, where $a_0,a_1,\ldots,a_n(\neq0)$ are complex
  constants and $t$ is the number of the distinct zeros of $P(z)$.
  Recently,
  Luo and Lin investigated more generally difference products of entire function and obtained the following result.
\medskip

{\noindent\bf Theorem D.}\cite[Theorem 1]{luolin}
  Let $f$ be a transcendental entire function of finite order.
For $n>t$, then $P(f)f(z+c)-a(z)$ has infinitely many zeros.

\medskip

Firstly, we give the following remark to show that the condition
$n>t$ in Theorem D is indispensable
 which is not given in \cite{luolin}.

\medskip

{\noindent\bf Remark.}  If $n=t=1$, Theorem D is not true, which can
be seen by the function $f(z)=e^{z}+1$, $e^{c}=-1$, hence
$f(z)f(z+c)-1=-e^{2z}$ has no zeros.

If $n=t=2$, Theorem D also is not true, which can be seen by
function $f(z)=\frac{1}{e^{z}}+1$, $e^{c}=-1$,
$P(z)=(z+\frac{-1+\sqrt{3}i}{2})(z+\frac{-1-\sqrt{3}i}{2})$, thus,
$P(f)f(z+c)-1=\frac{-1}{e^{3z}}$ has no zeros.

In fact, for any natural number $n=t$, we can construct an
counterexample to show Theorem D is not true by function
$f(z)=\frac{1}{e^{z}}+1$, $e^{c}=-1$,
$P(z)=(z-1-\frac{1}{d_{1}})\cdots(z-1-\frac{1}{d_{n-1}})$, where
$d_{i}\neq1, i=1,2,\ldots n-1$ are the distinct zero of $z^{n}-1=0$,
thus, we get $P(f)f(z+c)-1=\frac{-1}{e^{nz}}$ has no zeros.

\medskip

As the improvement of Theorem B, it is interesting to investigate
the zeros of derivatives of difference polynomials.
 The present
authors \cite[Theorem 1.1, Theorem 1.3]{liuliucao} have considered
the zeros of $[f^{n}f(z+c)]^{(k)}$ and $[f^{n}\Delta_{c}f]^{(k)}$,
the results can be stated as follows.

\medskip

{\noindent\bf Theorem E.} Let $f$ be a transcendental entire
function of finite order. If $n\geq k+2$, then
$[f(z)^nf(z+c)]^{(k)}-a(z)$ has infinitely many zeros. If $n\geq
k+3$, then $[f(z)^n\Delta_{c}f]^{(k)}-a(z)$ has infinitely many
zeros, unless $f$ is a periodic function with period $c$.

\medskip

In this paper, we continue to investigate the zeros of derivatives
of difference polynomials with more general forms and obtain the
following results as the improvements of the Theorem D and Theorem
E.

\begin{theorem}\label{june00}
  Let $f$ be a transcendental entire function of $\rho_{2}(f)<1$.
For $n\geq t(k+1)+1$, then $[P(f)f(z+c)]^{(k)}-a(z)$ has infinitely
many zeros.
\end{theorem}

{\noindent\bf Remark.} {\bf (1)}. Theorem \ref{june00} is an
improvement of Theorem E of the case $t=1$ and an improvement of
Theorem D of the case $k=0$.

{\bf (2)}. Theorem \ref{june00} is not valid for entire function
with $\rho_{2}(f)=1$, which can be seen by $f(z)=e^{e^{z}}$,
$P(z)=z^{n}$, $k\geq1$, $e^{c}=-n$, $a(z)$ is a nonconstant
polynomial, thus $[P(f)f(z+c)]^{(k)}-a(z)=-a(z)$ has finitely many
zeros.

{\bf (3)}. The condition of $a(z)\neq0$ can not be removed, which
can be seen by function $f(z)=e^{z}$, $P(z)=z^{n}$, $e^{c}=-1$, then
$[P(f)f(z+c)]^{(k)}=-(n+1)^{k}e^{(n+1)z}$ has no zeros.

\begin{theorem}\label{2010bb} Let $f$ be a transcendental entire function of $\rho_{2}(f)<1$, not a periodic function with period $c$.
 If $n\geq (t+1)(k+1)+1$, then
$[f(z)^n(\Delta_{c}f)^{s}]^{(k)}-a(z)$ has infinitely many zeros.
\end{theorem}

{\noindent\bf Remark.} The condition of $a(z)\neq0$ can not be
removed in Theorem \ref{2010bb}, which can be seen by function
$f(z)=e^{z}$, $P(z)=z^{n}$ $e^{c}=2$, then
$[P(f)\Delta_{c}f]^{(k)}=(n+1)^{k}e^{(n+1)z}$ has no zeros.

\medskip

For the case of transcendental meromorphic functions of Theorem
\ref{june00} and Theorem \ref{2010bb}, we obtian the next results.

\begin{theorem}\label{june2}
  Let $f$ be a transcendental meromorphic function of $\rho_{2}(f)<1$.
For $n\geq t(k+1)+5$, then $[P(f)f(z+c)]^{(k)}-a(z)$ has infinitely
many zeros.
\end{theorem}

{\noindent\bf Remark.} Theorem \ref{june2} also partly answer the
question raised by Luo and Lin \cite[P. 448]{luolin}.

\begin{theorem}\label{june3}
  Let $f$ be a transcendental meromorphic function of
  $\rho_{2}(f)<1$.
For $n\geq(t+2)(k+1)+3+s$, then $[P(f)(\Delta_{c}f)^{s}]^{(k)}-a(z)$
has infinitely many zeros.
\end{theorem}

\begin{corollary}\label{co2} Let $P(z), Q(z), H(z)$ be
nonzero polynomials.  If $H(z)$ is a nonconstant polynomial, then
the nonlinear difference-differential equation
\begin{equation}\label{00011110}
[P(f)f(z+c)]^{(k)}-P(z)=Q(z)e^{H(z)}
\end{equation}
has no transcendental entire $($meromorphic$)$ solution of
$\rho_{2}(f)<1$, provided that $n\geq t(k+1)+1$ $(n\geq t(k+1)+5)$.
If $H(z)$ is a constant, then (\ref{00011110}) has no transcendental
entire solutions of $\rho_{2}(f)<1$, and (\ref{00011110}) has no
transcendental meromorphic solutions of $\rho_{2}(f)<1$ provided
that $n\geq2$.
\end{corollary}

\begin{corollary}\label{cococo} Let $P(z), Q(z), H(z)$ be
nonzero polynomials. If $H(z)$ is a nonconstant polynomial, then the
nonlinear difference-differential equation
\begin{equation}\label{00022220}
[P(f)(\Delta_{c}f)^{s}]^{(k)}-P(z)=Q(z)e^{H(z)}
\end{equation}
has no transcendental entire $($meromorphic$)$solution of
$\rho_{2}(f)<1$, provided that $n\geq (t+1)(k+1)+s+1$
$(n\geq(t+2)(k+1)+3+s)$. If $H(z)$ is a constant, then
(\ref{00022220}) has no transcendental entire solutions of
$\rho_{2}(f)<1$, and (\ref{00022220}) has no transcendental
meromorphic solutions of  $\rho_{2}(f)<1$ provided that $n\geq3$,
unless $f$ is a periodic function with period $c$.
\end{corollary}

About the uniqueness of difference products of entire functions,
some results can be found in \cite{liucao, liuliucao, liuzhangyang,
luolin, Qi, zhang11}. The main purpose is to obtain the
relationships between $f$ and $g$ when $P(f)f(z+c)$ and $P(g)g(z+c)$
sharing one common value. In fact, two special types $P(z)=z^{n}$
and $P(z)=z^{n}(z^{m}-1)$ always be considered. Luo and Lin
\cite[Theorem 2]{luolin} also considered the general case of $P(z)$.
Corresponding to the above theorems of this paper, it is necessary
to consider the uniqueness of derivative of difference polynomials
sharing one common value. The present authors \cite[Theorem
1.5]{liuliucao} have considered the uniqueness about
$[f^{n}f(z+c)]^{(k)}$ and $[g^{n}g(z+c)]^{(k)}$ sharing one common
value, the result can be stated as follows.

\medskip

{\noindent\bf Theorem F.} Let $f(z)$ and $g(z)$ be transcendental
entire functions of finite order, $n\geq 2k+6$. If
$[f(z)^nf(z+c)]^{(k)}$ and $[g(z)^ng(z+c)]^{(k)}$ share the value 1
CM, then either $f(z)=c_{1}e^{Cz}$, $g(z)=c_{2}e^{-Cz}$, where
$c_{1}, c_{2}$ and $C$ are constants satisfying
$(-1)^{k}(c_{1}c_{2})^{n}[(n+1)C]^{2k}=1$ or $f=tg$, where
$t^{n+1}=1$.

\medskip

In this paper, we consider the entire functions of $\rho_{2}(f)<1$
and get the following theorems.

\begin{theorem}\label{2010coro}Let $f(z)$ and $g(z)$ be transcendental entire functions
of $\rho_{2}(f)<1$, $n\geq 2k+m+6$. If
$[f^{n}(f^{m}-1)f(z+c)]^{(k)}$ and $[g^{n}(g^{m}-1)g(z+c)]^{(k)}$
share the value 1 CM, then $f=tg$, where $t^{n+1}=t^{m}=1$.
\end{theorem}

\begin{theorem} \label{2010d} The conclusion of Theorem \ref{2010coro} is also valid, if $n\geq
5k+4m+12$ and $[f^{n}(f^{m}-1)f(z+c)]^{(k)}$ and
$[g^{n}(g^{m}-1)g(z+c)]^{(k)}$ share the value 1 IM.
\end{theorem}

\section{Some Lemmas}
For a finite order transcendental meromorphic function $f$, the
difference logarithmic derivative lemma, given by Chiang and Feng
\cite[Corollary 2.5]{chiang}, Halburd and Korhonen \cite[Theorem
2.1]{R.G}, \cite[Theoem 5.6]{hal}, plays an important part in
considering the difference Nevanlinna theory. Afterwards, R. G.
Halburd, R. J. Korhonen and K. Tohge improved the condition of
growth from finite order to $\rho_{2}(f)<1$ as follows.

\begin{lemma}\label{lem2.1}\cite[Theorem 5.1]{tohge}
Let $f$ be a transcendental meromorphic function of $\rho_{2}(f)<1$,
$\varsigma<1$, $\varepsilon$ is a enough small number. Then
\begin{equation}\label{2.1}
m\left(r,\frac{f(z+c)}{f(z)}\right)=o\left(\frac{T(r,f)}{r^{1-\varsigma-\varepsilon}}\right)=S(r,
f),
\end{equation}
for all $r$ outside of a set of finite logarithmic measure.
 \end{lemma}

\begin{lemma}\label{lem2.2}\cite[Lemma 8.3]{tohge}
Let $T:[0,+\infty)\rightarrow[0,+\infty)$ be a non-decreasing
continuous function and let $s\in(0,\infty).$ If the hyper order of
$T$ is strictly less that one, i.e.,
\begin{eqnarray}\label{aiyundong}
\limsup_{r\rightarrow\infty}\frac{\log\log T(r)}{\log
r}=\varsigma<1,
\end{eqnarray}
and $\delta\in(0,1-\varsigma),$ then
\begin{equation}\label{ai}
T(r+s)=T(r)+o\left(\frac{T(r)}{r^{\delta}}\right)
\end{equation}
for all $r$ runs to infinity outside of a set of finite logarithmic
measure.
\end{lemma}

Thus, from Lemma \ref{lem2.2}, we get the following lemma.
\begin{lemma}\label{li}
Let $f(z)$ be a transcendental meromorphic function of
$\rho_{2}(f)<1$. Then,
\begin{equation}\label{tiancai}
T(r,f(z+c))=T(r,f)+S(r,f)
\end{equation}
and
\begin{equation}\label{tiancai}
N(r,f(z+c))=N(r,f)+S(r,f), \qquad
N\left(r,\frac{1}{f(z+c)}\right)=N(r,\frac{1}{f})+S(r,f).
\end{equation}
\end{lemma}

Combining the method of proof of \cite[Lemma 5]{luolin} with Lemma
\ref{lem2.1}, we can get the following Lemma
\ref{julyshang}--Lemma\ref{duan}.
\begin{lemma}\label{julyshang}
Let $f(z)$ be a transcendental entire function of $\rho_{2}(f)<1$.
If $F=P(f)f(z+c)$, then \begin{equation}\label{liu1}
T(r,F)=T(r,P(f)f(z))+S(r,f)=(n+1)T(r,f)+S(r,f).
\end{equation}
\end{lemma}

\begin{lemma}\label{june4444}
Let $f(z)$ be a transcendental meromorphic function of
$\rho_{2}(f)<1$. If $F=P(f)f(z+c)$, then
\begin{equation}\label{liu1}
(n-1)T(r,f)+S(r,f)\leq T(r,F)\leq (n+1)T(r,f)+S(r,f).
\end{equation}
\end{lemma}

\begin{proof}
Since $F(z)=P(f)f(z+c)$, then
\begin{equation}\label{51}
\frac{1}{P(f)f}=\frac{1}{F}\frac{f(z+c)}{f(z)}.
\end{equation}
Using the first and second main theorem, Lemma \ref{lem2.1} and the
standard Valrion-Monko's theorem \cite{Mohon}, from (\ref{51}), we
get
\begin{eqnarray}\label{511}
(n+1)T(r,f)&\leq& T(r,F(z))+T(r,\frac{f(z+c)}{f(z)})+O(1)\nonumber\\
&\leq& T(r,F(z))+m(r,\frac{f(z+c)}{f(z)})+N(r,\frac{f(z+c)}{f(z)})+O(1)\nonumber\\
&\leq&T(r,F(z))+N(r,\frac{f(z+c)}{f(z)})+S(r,f)\nonumber\\
&\leq&T(r,F(z))+2T(r,f)+S(r,f),
\end{eqnarray}
hence, we get $T(r,F)\geq(n-1)T(r,f)+S(r,f)$. It is easy to get
$T(r,F)\leq (n+1)T(r,f)+S(r,f)$. Thus, (\ref{liu1}) follows.
\end{proof}

{\noindent\bf Remark.} The inequality (\ref{liu1}) can not be
improved by the following two examples. If $f(z)=\tan z$,
$P(z)=z^{n}$, $c_{1}=\frac{\pi}{2}$, then
$$T(r,P(z)f(z+c_{1}))=-\tan^{n-1} z=(n-1)T(r,f)+S(r,f).$$
If $f(z)=\tan z$, $P(z)=z^{n}$, $c_{2}=\pi$, then
$$T(r,P(z)f(z+c_{2}))=\tan^{n+1}z=(n+1)T(r,f)+S(r,f).$$

\begin{lemma}\label{li222li}
Let $f(z)$ be a transcendental entire function of $\rho_{2}(f)<1$.
Then,
\begin{equation}\label{rho}
nT(r,f)+S(r,f)\leq T(r,P(f)[f(z+c)-f(z)]^{s})\leq
(n+s)T(r,f)+S(r,f).
\end{equation}
\end{lemma}

{\noindent\bf Remark.} The inequality (\ref{rho}) can not be
improved by the following two examples. If $f(z)=e^{z}$, $e^{c}=2$,
then
$$T(r,f(z)^{n}[f(z+c)-f(z)]^{s})=T(r,e^{(n+s)z})=(n+s)T(r,f)+S(r,f).$$
If $f(z)=e^{z}+z$, $c=2\pi i$, then
$$T(r,f(z)^{n}[f(z+c)-f(z)^{s}])=T(r,(2\pi i)^{s}[e^{z}+z]^{n})=nT(r,f)+S(r,f).$$

\begin{lemma}\label{duan}
Let $f(z)$ be a transcendental meromorphic function of
$\rho_{2}(f)<1$. Then,
\begin{equation}\label{rho0000}
(n-s)T(r,f)+S(r,f)\leq
T(r,P(f)[f(z+c)-f(z)]^{s})\leq(n+2s)T(r,f)+S(r,f).
\end{equation}
\end{lemma}

The following lemma is needed for the proof of Theorem
\ref{2010coro}. For the case of $k=0$, $m=1$, $f$ and $g$ are
transcendental entire functions of finite order, the proof can be
found in \cite[The proof of Theorem 6]{zhang11}.

\begin{lemma}\label{nan}
Let $f$ and $g$ be transcendental entire functions of
$\rho_{2}(f)<1$, and $c$ be a nonzero constant. If $n\geq m+5$ and
\begin{equation}\label{haha}
[f^{n}(f^{m}-1)f(z+c)]^{(k)}=[g^{n}(g^{m}-1)g(z+c)]^{(k)},
\end{equation} then $f=tg$, and $t^{n+1}=t^{m}=1$.
\end{lemma}

\begin{proof} From (\ref{haha}), we get $f^{n}(f^{m}-1)f(z+c)=g^{n}(g^{m}-1)g(z+c)+Q(z)$, where $Q(z)$ is a
polynomial of degree at most $k-1.$ If $Q(z)\not\equiv0,$ then we
have
$$\frac{f^{n}(f^{m}-1)f(z+c)}{Q(z)}=\frac{g^{n}(g^{m}-1)g(z+c)}{Q(z)}+1.$$
From the second main theorem of Nevanlinna and Lemma
\ref{julyshang}, we have
\begin{eqnarray}\label{bbb}
(n+m+1)T(r,f)&=&T(r,\frac{f^{n}(f^{m}-1)f(z+c)}{Q(z)})+S(r,f)\nonumber\\&\leq&
\overline{N}(r,\frac{f^{n}(f^{m}-1)f(z+c)}{Q(z)})+\overline{N}(r,\frac{Q(z)}{f^{n}(f^{m}-1)f(z+c)})\nonumber\\&&+\overline{N}(r,\frac{Q(z)}{g^{n}(g^{m}-1)g(z+c)})+S(r,f)
\nonumber\\&\leq&\overline{N}(r,\frac{1}{f^{n}(f^{m}-1)})+\overline{N}(r,\frac{1}{f(z+c)})+\overline{N}(r,\frac{1}{g^{n}(g^{m}-1)})
\nonumber\\&&+\overline{N}(r,\frac{1}{g(z+c)})+S(r,f)
\nonumber\\&\leq&(m+2)T(r,f)+(m+2)T(r,g)+S(r,f)+S(r,g).
\end{eqnarray}
Similarly as above, we have
\begin{eqnarray*}
(n+m+1)T(r,g)\leq(m+2)T(r,f)+(m+2)T(r,g)+S(r,f)+S(r,g).
\end{eqnarray*}
Thus, we get
\begin{eqnarray*}
(n+m+1)[T(r,f)+T(r,g)]\leq2(m+2)[T(r,f)+T(r,g)]+S(r,f)+S(r,g).
\end{eqnarray*}
which is a contradiction with $n\geq m+5$. Hence, we get
$P(z)\equiv0$, which implies that
\begin{equation}\label{buhao}
f^{n}(f^{m}-1)f(z+c)=g^{n}(g^{m}-1)g(z+c).
\end{equation}
Let $G(z)=\frac{f(z)}{g(z)}.$ Assume that $G(z)$ is nonconstant.
From (\ref{buhao}), we have
      \begin{eqnarray}\label{6.1}
        g(z)^{m}=\frac{G(z)^nG(z+c)-1}{G(z)^{n+m}G(z+c)-1}.
           \end{eqnarray}
If 1 is a Picard value of $G(z)^{n+m}G(z+c)$, then applying the
second main theorem, we get
\begin{eqnarray}\label{6.2}
T(r, G^{n+m}G(z+c))&\leq&\overline{N}(r, G^{n+m}G(z+c))+\overline{N}(r,\frac{1}{G^{n+m}G(z+c)})\nonumber\\
&+& \overline{N}(r, \frac{1}{G^{n+m}G(z+c)-1})+S(r,G)\nonumber\\
&\leq& 2T(r, G(z))+2T(r, G(z+c))+S(r,G)\nonumber\\
&\leq& 4T(r, G(z))+S(r,G).
\end{eqnarray}
Combining (\ref{6.2}) with Lemma \ref{june4444}, we have
$(n+m-1)T(r,G)\leq4T(r, G(z))+S(r,G),$ which is a contradiction with
$n\geq m+5$. Therefore, 1 is not a Picard value of
$G(z)^{n+m}G(z+c)$. Thus, there exists $z_0$ such that
$G(z_0)^{n+m}G(z_0+c)=1$. The following, we may distinguish two
cases.

{\bf Case 1.} $G(z)^{n+m}G(z+c)\not\equiv1$. From (\ref{6.1}) and
$g(z)$ is an entire function, then we get $G(z_0)^nG(z_0 + c)=1$,
thus $G(z_0)^{m}=1$. Therefore,
\begin{eqnarray}\label{6.4}
\overline{N}(r, \frac{1}{G^{n+m}G(z+c)-1})\leq \overline{N}(r,
\frac{1}{G^{m}-1})\leq m T(r, G)+ S(r, G).
\end{eqnarray}
By (\ref{6.4}) and Lemma \ref{li}, applying the second main theorem,
we get
\begin{eqnarray}\label{6.5}
\begin{split}
T(r, G^{n+m}G(z+c))&\leq \overline{N}(r, G^{n+m}G(z+c))+\overline{N}(r, \frac{1}{G^{n+m}G(z+c)})\\
&+ \overline{N}(r, \frac{1}{G^{n+m}G(z+c)-1})+S(r, G)\\
&\leq (m+2)T(r,G(z))+2T(r, G(z+c))+S(r, G)\\
& \leq (m+4)T(r, G(z))+S(r, G).
\end{split}
\end{eqnarray}
On the other hand, we have
\begin{eqnarray}\label{6.6}
(n+m)T(r, G)&=&T(r, G^{n+m})\nonumber\\&\leq&T(r, G^{n+m}G(z+c))+T(r, G(z+c))+O(1)\nonumber\\
&\leq& (m+5)T(r, G(z))+S(r, G),
\end{eqnarray}
which contradicts $n\geq m+5\geq6$.

{\bf Case 2.} $G(z)^{n+m}G(z+c)\equiv1$, thus,
\begin{eqnarray}\label{6.6}
(n+m)T(r,G)&=&T(r,G(z+c))+S(r,G)\nonumber\\&=&T(r, G(z))+S(r, G),
\end{eqnarray}
 which also is a contradiction with $n\geq m+5$.
Thus, $G$ must be a constant and $f(z)=tg(z)$, where $t$ is a
non-zero constant. From $f^{n}(f^{m}-1)f(z+c)\equiv
g^{n}(g^{m}-1)g(z+c)$, we know that $t^{m}=1$ and $t^{n+1}=1$, $n,
m$ are positive integers.
\end{proof}

The following result is related to the growth of solutions of linear
difference equation and is needed for the proof of Lemma
\ref{jintian}, was given by Li and Gao \cite[Theorem 2.1]{lisheng}.
Here, we give the version with small changes of the type of equation
(\ref{wang}), the proof are similar.

\begin{lemma}\label{ligao}
Let $a_{0}(z), a_{1}(z),\cdots,a_{n}(z), b(z)$ be polynomials such
that $a_{0}(z)a_{n}(z)\not\equiv0$, let $c_{j}$ be constants and
$$\deg(\sum_{\deg a_{j}=d} a_{j})=d,$$ where $d=\max_{0\leq j\leq
n}\{\deg a_{j}\}$. If $f(z)$ is a transcendental meromorphic
solution of
\begin{equation}\label{wang}
\sum_{j=0}^{n}a_{j}(z)f(z+c_{j})=b(z),
\end{equation} then $\rho(f)\geq1$.
\end{lemma}

\begin{lemma}\label{jintian}
If $n\geq k+1$, then there are no transcendental entire functions
$f$ and $g$ with hyper order less than one, satisfying
\begin{equation}\label{772227}
[f^{n}(f^{m}-1)f(z+c)]^{(k)}\cdot [g^{n}(g^{m}-1)g(z+c)]^{(k)}=1.
\end{equation}
\end{lemma}

\begin{proof} Assume that $f$ and $g$ satisfy (\ref{772227}) and $f$ and $g$ are
transcendental entire functions of hyper order less than one. Since
$n\geq k+1$, from (\ref{772227}), we get $f$ and $g$ have no zeros.
Thus, $f(z)=e^{b(z)}$ and $g(z)=e^{d(z)}$, where $b(z), d(z)$ are
entire functions with order less than one. Thus, substitute $f$ and
$g$ into (\ref{772227}), we get
$$[e^{nb(z)}(e^{mb(z)}-1)e^{b(z+c)}]^{(k)}[e^{nd(z)}(e^{md(z)}-1)e^{d(z+c)}]^{(k)}=1$$
Let $(n+m)b(z)+b(z+c)=B_{1}(z)$, $nb(z)+b(z+c)=B_{2}(z)$ and
$(n+m)d(z)+d(z+c)=D_{1}(z)$, $nd(z)+d(z+c)=D_{1}(z)$.
\medskip

If $k=1$, we have
$$[B'_{1}(z)e^{B_{1}(z)}-B'_{2}(z)e^{B_{2}(z)}][D'_{1}(z)e^{D_{1}(z)}-D'_{2}(z)e^{D_{2}(z)}]=1,$$
which implies that
$e^{B_{2}(z)}[B'_{1}(z)e^{B_{1}(z)-B_{2}(z)}-B'_{2}(z)]$ has no
zeros. If $B'_{1}\not=0$, remark that 0 is the Picard exceptional
value of $e^{B_{1}(z)-B_{2}(z)}$, then we get $B'_{2}(z)$ must be
zero, thus $B_{2}$ must be a constant. From Lemma \ref{ligao} and
$nb(z)+b(z+c)=B_{2}$, we get $\rho(b(z))\geq1$, thus
$\rho_{2}(f)\geq1$, which is a contradiction. If $B'_{1}=0$, then
$B_{1}$ must be a constant, which also induces that
$\rho_{2}(f)\geq1$, a contradiction.
\medskip

If $k=2$, by calculation, then we have
$e^{B_{2}(z)}[(B''_{1}(z)+B'^{2}_{1}(z))e^{B_{1}(z)-B_{2}(z)}-(B''_{2}(z)+B'^{2}_{2}(z))]$
has no zeros. If $B''_{1}+B'^{2}_{1}\not=0$, then
$B''_{2}+B'^{2}_{2}=0$. If $B_{2}$ is transcendental entire, then we
get
$$m(r,B'_{2})=m(r,\frac{B''_{2}}{B'_{2}})=S(r,B'_{2}),$$ which is a
contradiction with $B'_{2}$ is transcendental entire. If  $B_{2}$ is
a polynomial, from Lemma \ref{ligao}, which also induces that
$\rho_{2}(f)\geq1$, a contradiction. If $B''_{1}+B'^{2}_{1}=0$,
similar as above, we get a contradiction. For any $k\geq2$, using
the similar method as above, we can get the proof of Lemma
\ref{jintian}.
\end{proof}

Let $p$ be a positive integer and $a\in\mathbb{C}$. We denote by
$N_{p}(r,\frac{1}{f-a})$ the counting function of the zeros of $f-a$
where an $m$-fold zero is counted $m$ times if $m\leq p$ and $p$
times if $m>p$.

\begin{lemma}\label{jintianmingtian}
Let $f$ be a nonconstant meromorphic function, and $p, k$ be
positive integers. Then
\begin{equation}\label{eeed}
T(r,f^{(k)})\leq T(r,f)+k\overline{N}(r,f)+S(r,f).
\end{equation}
\begin{eqnarray}\label{yu}
N_{p}(r, \frac{1}{f^{(k)}})\leq T(r,f^{(k)})-T(r,f) + N_{p+k}(r,
\frac{1}{f})+S(r, f),
\end{eqnarray}
\begin{eqnarray}\label{yuyu}
N_{p}(r,\frac{1}{f^{(k)}})\leq k \overline{N}(r,f)+
N_{p+k}(\frac{1}{f})+S(r, f),
\end{eqnarray}
\end{lemma}

\begin{lemma}\label{55}\cite[Lemma 3]{hua} Let $F$ and $G$ be nonconstant
meromorphic functions. If $F$ and $G$ share 1 CM, then one of the
following three cases holds:
\begin{itemize}
\item[(i)]$\max\{T(r,F ),T(r,G)\}\leq N_{2}(r,
\frac{1}{F})+N_{2}(r,F)+N_{2}(r,\frac{1}{G})+N_{2}(r,G)+S(r,F)+S(r,G),$
\item[(ii)] $F=G$,
\item[(iii)]$F\cdot G=1$.
\end{itemize}
\end{lemma}
For the proof of Theorem \ref{2010d}, we need the following lemma.
\begin{lemma}\label{66}\cite[Lemma 2.3]{xu} Let $F$ and $G$ be nonconstant
meromorphic functions sharing the value 1 IM. Let
$$H=\frac{F''}{F'}-2\frac{F'}{F-1}-\frac{G''}{G'}+2\frac{G'}{G-1}.$$
If $H\not\equiv0$, then
\begin{eqnarray}\label{h}
T(r,F)+T(r,G)&\leq& 2\left(N_{2}(r,
\frac{1}{F})+N_{2}(r,F)+N_{2}(r,\frac{1}{G})+N_{2}(r,G)\right)\nonumber\\&+&
3\left(\overline{N}(r,F)+\overline{N}(r,\frac{1}{F})+\overline{N}(r,G)+\overline{N}(r,\frac{1}{G})\right)\nonumber\\&+&S(r,F)+S(r,G).
\end{eqnarray}
\end{lemma}

\section{Proofs of Theorem \ref{june00} and Theorem \ref{2010bb}}
Let $F(z)=P(f)f(z+c).$  From Lemma \ref{julyshang}, we know that
$F(z)$ is not a constant, and $S(r,F)=S(r,F^{(k)})=S(r,f)$ follows.
Assume that $F(z)^{(k)}-\alpha(z)$ has only finitely many zeros,
combining the second main theorem for three small functions
\cite[Theorem 2.5]{hayman3} and (\ref{yu}) with $f$ is a
transcendental entire function, then we get
\begin{eqnarray}\label{9324}
T(r,F^{(k)})&\leq&
\overline{N}(r,F^{(k)})+\overline{N}(r,\frac{1}{F^{(k)}})+\overline{N}(r,\frac{1}{F^{(k)}-\alpha(z)})+S(r,F^{(k)})\nonumber\\&&
\leq
N_{1}(r,\frac{1}{F^{(k)}})+\overline{N}(r,\frac{1}{F^{(k)}-\alpha(z)})+S(r,F^{(k)})\nonumber\\&&
\leq T(r,F^{(k)})-T(r,F) + N_{k+1}(r, \frac{1}{F})+S(r,F^{(k)}).
\end{eqnarray}
Combining (\ref{liu1}) with (\ref{9324}), it implies that
\begin{eqnarray}\label{94}
(n+1)T(r,f)+S(r,f)&=&T(r,F)\leq N_{k+1}(r, \frac{1}{F})+S(r,
f)\nonumber\\&\leq&t(k+1)\overline{N}(r,\frac{1}{f})+N(r,\frac{1}{f(z+c)})+S(r,
f)\nonumber\\&\leq&[t(k+1)+1]T(r,f)+S(r, f),
\end{eqnarray}
which is a contradiction with $n\geq t(k+1)+1$. Thus, Theorem
\ref{june00} is proved. Set $G(z)=P(f)[\Delta_{c}f]^{s}.$ If
$G(z)^{(k)}-\alpha(z)$ has only finitely many zeros, using the
similar method as above, from Lemma \ref{li222li}, then we get
\begin{eqnarray}\label{949999}
n T(r,f)+S(r,f)&\leq&T(r,G)\leq N_{k+1}(r, \frac{1}{G})+S(r,
f)\nonumber\\&\leq&t(k+1)\overline{N}(r,\frac{1}{f})+(k+1)\overline{N}(r,\frac{1}{f(z+c)-f(z)})+S(r,
f)\nonumber\\&\leq&(t+1)(k+1)T(r,f)+S(r, f),
\end{eqnarray}
which is a contradiction with $n\geq(t+1)(k+1)+1$.
 Thus, we get the
proof of Theorem \ref{2010bb}.

\section{Proofs of Theorem \ref{june2} and Theorem \ref{june3}}
Let $F(z)=P(f)f(z+c).$  From Lemma \ref{june4444}, we know that
$F(z)$ is not a constant, and $S(r,F)=S(r,F^{(k)})=S(r,f)$ follows.
Assume that $F(z)^{(k)}-\alpha(z)$ has only finitely many zeros,
combining the second main theorem for three small functions
\cite[Theorem 2.5]{hayman3} and (\ref{yu}) with $f$ is a
transcendental entire function, then we get
\begin{eqnarray}\label{93}
T(r,F^{(k)})&\leq&
\overline{N}(r,F^{(k)})+\overline{N}(r,\frac{1}{F^{(k)}})+\overline{N}(r,\frac{1}{F^{(k)}-\alpha(z)})+S(r,F^{(k)})\nonumber\\&
\leq&
\overline{N}(r,f)+\overline{N}(r,f(z+c))+N_{1}(r,\frac{1}{F^{(k)}})+\overline{N}(r,\frac{1}{F^{(k)}-\alpha(z)})+S(r,F^{(k)})\nonumber\\&
\leq& 2T(r,f)+T(r,F^{(k)})-T(r,F) + N_{k+1}(r,
\frac{1}{F})+S(r,F^{(k)}).
\end{eqnarray}
Combining (\ref{liu1}) with (\ref{93}), it implies that
\begin{eqnarray}\label{94}
(n-1)T(r,f)+S(r,f)&\leq&T(r,F)\leq 2T(r,f)+N_{k+1}(r,
\frac{1}{F})+S(r,
f)\nonumber\\&\leq&t(k+1)\overline{N}(r,\frac{1}{f})+N(r,\frac{1}{f(z+c)})+2T(r,f)+S(r,
f)\nonumber\\&\leq&[t(k+1)+3]T(r,f)+S(r, f),
\end{eqnarray}
which is a contradiction with $n\geq t(k+1)+5$. Thus, Theorem
\ref{june2} is proved. Set $G(z)=P(f)[\Delta_{c}f]^{s}.$ If
$G(z)^{(k)}-\alpha(z)$ has only finitely many zeros, using the
similar method as above, from Lemma \ref{li222li}, then we get
\begin{eqnarray}\label{949999}
(n&-&s)T(r,f)+S(r,f)\leq T(r,G)\leq 2T(r,f)+N_{k+1}(r,
\frac{1}{G})+S(r,
f)\nonumber\\&\leq&2T(r,f)+t(k+1)\overline{N}(r,\frac{1}{f})+(k+1)\overline{N}(r,\frac{1}{f(z+c)-f(z)})+S(r,
f)\nonumber\\&\leq&[(t+2)(k+1)+2]T(r,f)+S(r, f),
\end{eqnarray}
which is a contradiction with $n\geq(t+2)(k+1)+3+s$.
 Thus, we get the
proof of Theorem \ref{june3}.

\section{Proof of Theorem \ref{2010coro}}
Let $F=[f^{n}(f^{m}-1)f(z+c)]^{(k)}$,
$G=[g^{n}(g^{m}-1)g(z+c)]^{(k)}$. Thus $F$ and $G$ share the value 1
CM. From (\ref{eeed}) and $f$ is a transcendental entire function,
then
\begin{eqnarray}\label{0000}
T(r,F)\leq T(r,f^{n}(f^{m}-1)f(z+c))+S(r,P(f)f(z+c)).
\end{eqnarray}
 Combining
(\ref{0000}) with (\ref{julyshang}), we have $S(r,F)=S(r,f)$. We
also have $S(r,G)=S(r,g)$ from the same reason as above. From
(\ref{yu}), we obtain
\begin{eqnarray}\label{1}
N_{2}(r,\frac{1}{F})&=&
N_{2}\left(r,\frac{1}{[f^{n}(f^{m}-1)f(z+c)]^{(k)}}\right)\nonumber\\&\leq&
T(r,F)-T(r,f^{n}(f^{m}-1)f(z+c))\nonumber\\&+&N_{k+2}(r,\frac{1}{f^{n}(f^{m}-1)f(z+c)})+S(r,f).
\end{eqnarray}
Thus, from Lemma \ref{julyshang} and (\ref{1}), we get
\begin{eqnarray}\label{3}
(n+m+1)T(r,f)&=&T(r,f^{n}(f^{m}-1)f(z+c))+S(r,f)\nonumber\\&\leq&
T(r,F)-N_{2}(r,\frac{1}{F})+N_{k+2}(r,\frac{1}{f^{n}(f^{m}-1)f(z+c)})+S(r,f).
\end{eqnarray}
From (\ref{yuyu}), we obtain
\begin{eqnarray}\label{111}
N_{2}(r,\frac{1}{F})&\leq&
N_{k+2}(r,\frac{1}{f^{n}(f^{m}-1)f(z+c)})+S(r,f) \nonumber\\&\leq&
(k+2)
N(r,\frac{1}{f})+N(r,\frac{1}{f^{m}-1})+N(r,\frac{1}{f(z+c)})+S(r,f)\nonumber\\&\leq&
(k+m+3)T(r,f)+S(r,f).
\end{eqnarray}
Similarly as above, we obtain
\begin{eqnarray}\label{33}
(n+m+1)T(r,g)&\leq&
T(r,G)-N_{2}(r,\frac{1}{G})\nonumber\\&+&N_{k+2}(r,\frac{1}{g^{n}(g^{m}-1)g(z+c)})+S(r,g).
\end{eqnarray}
and
\begin{eqnarray}\label{2}
N_{2}(r,\frac{1}{G})\leq(k+m+3)T(r,g)+S(r,g).
\end{eqnarray}
If the (i) of Lemma \ref{55} is satisfied, implies that
$$\max\{T(r,F),T(r,G)\}\leq
N_{2}(r,\frac{1}{F})+N_{2}(r,\frac{1}{G})+S(r,F)+S(r,G).$$ Thus,
combining above with (\ref{3})--(\ref{2}), we obtain
\begin{eqnarray}\label{66}
(n+m+1)[T(r,f)+T(r,g)]&\leq&
2N_{k+2}(r,\frac{1}{f^{n}(f^{m}-1)f(z+c)})\nonumber\\&+&2N_{k+2}(r,\frac{1}{g^{n}(g^{m}-1)g(z+c)})+S(r,f)+S(r,g)\nonumber\\&\leq&
2(k+m+3)[T(r,f)+T(r,g)]+S(r,f)+S(r,g),
\end{eqnarray}
which is a contradiction with  $n\geq 2k+m+6$. Hence, $F=G$ or
$F\cdot G=1$. From Lemma \ref{nan} and Lemma \ref{jintian}, we get
$f=tg$ for $t^{m}=t^{n+1}=1$. Thus, we get the proof of Theorem
\ref{2010coro}.

\section{Proof of Theorem \ref{2010d}}
Let $F=[f^{n}(f^{m}-1)f(z+c)]^{(k)}$,
$G=[g^{n}(g^{m}-1)g(z+c)]^{(k)}$. We will show that $F=G$ or $F\cdot
G=1$ under the conditions of Theorem \ref{2010d}. Assume that
$H\not\equiv0$, from (\ref{h}), we get
\begin{eqnarray}\label{hh}
T(r,F)+T(r,G)&\leq&2\left(N_{2}(r,
\frac{1}{F})+N_{2}(r,\frac{1}{G})\right)+
3\left(\overline{N}(r,\frac{1}{F})+\overline{N}(r,\frac{1}{G})\right)\nonumber\\&+&S(r,F)+S(r,G).
\end{eqnarray}
Combining above with (\ref{3})--(\ref{2}) and (\ref{yuyu}), we
obtain
\begin{eqnarray*}\label{hhgg}
(n&+&m+1)(T(r,f)+T(r,g))\leq T(r,F)+T(r,G)+N_{k+2}(r,\frac{1}{f^{n}(f^{m}-1)f(z+c)})\nonumber\\
&+&N_{k+2}(r,\frac{1}{g^{n}(g^{m}-1)g(z+c)})-N_{2}(r,\frac{1}{F})-N_{2}(r,\frac{1}{G})+S(r,f)+S(r,g)\nonumber\\&\leq&2N_{k+2}(r,\frac{1}{f^{n}(f^{m}-1)f(z+c)})+2N_{k+2}(r,\frac{1}{g^{n}(g^{m}-1)g(z+c)})
\nonumber\\&+&3\left(\overline{N}(r,\frac{1}{F})+\overline{N}(r,\frac{1}{G})\right)+S(r,f)+S(r,g)\nonumber\\&\leq&(5k+5m+12)[T(r,f)+T(r,g)]+S(r,f)+S(r,g),
\end{eqnarray*}
which is a contradiction with $n\geq 5k+4m+12$. Thus, we get
$H\equiv0$. The following proof is trivial, the original idea is
devoting to Yang and Yi \cite{yang and yi}. Here, we give the
complete proof. Integrating $H$ twice, we obtain
\begin{equation}\label{leilee}
F=\frac{(b+1)G+(a-b-1)}{bG+(a-b)}, G=\frac{(a-b-1)-(a-b)F}{Fb-(b+1)}
\end{equation}
which implies that $T(r,F)=T(r,G)+O(1)$.  We divide into three cases
as follows:

{\bf Case 1}. $b\not=0, -1$. If $a-b-1\not=0$, then by
(\ref{leilee}), we get
\begin{equation}\label{leilelelele}
\overline N(r,\frac{1}{F})=\overline
N\left(r,\frac{1}{G-\frac{a-b-1}{b+1}}\right).
\end{equation}
By the Nevanlinna second main theorem, (\ref{yu}) and (\ref{yuyu}),
we have
\begin{eqnarray}\label{ddd}
(n+m+1)T(r,g)&\leq&T(r,G)+N_{k}(r,\frac{1}{g^{n}(g^{m}-1)g(z+c)})-N(r,\frac{1}{G})+S(r,g)\nonumber\\&\leq&N_{k}(r,
\frac{1}{g^{n}(g^{m}-1)g(z+c)})+\overline{N}\left(r,
\frac{1}{G-\frac{a-b-1}{b+1}}\right)+S(r,g)\nonumber\\&\leq&
(k+m+1)T(r,g)+(k+m+2)T(r,f)+S(r,f)+S(r,g).
\end{eqnarray}
Similarly, we get \begin{eqnarray*}\label{ddddddd}
(n+m+1)T(r,f)\leq(k+m+1)T(r,f)+(k+m+2)T(r,g)+S(r,f)+S(r,g).
\end{eqnarray*}
Thus, from (\ref{ddd}) and above, then
$$(n+m+1)[T(r,f)+T(r,g)]\leq(2k+2m+3)[T(r,f)+T(r,g)]+S(r,f)+S(r,g),$$
which is a contradiction with $n\geq 5k+4m+12$. Thus, $a-b-1=0$,
then
\begin{equation}\label{leilele}
F=\frac{(b+1)G}{bG+1}.
\end{equation}
Since $F$ is an entire function and (\ref{leilele}), then
$\overline{N}(r, \frac{1}{G+\frac{1}{b}})=0$. Using the same method
as above, we get
\begin{eqnarray}\label{dddllll}
(n+m+1)T(r,g)&\leq&T(r,G)+N_{k}(r,\frac{1}{g^{n}(g^{m}-1)g(z+c)})-N(r,\frac{1}{G})+S(r,g)\nonumber\\&\leq&N_{k}(r,
\frac{1}{g^{n}(g^{m}-1)g(z+c)})+\overline{N}\left(r,
\frac{1}{G+\frac{1}{b}}\right)+S(r,g)\nonumber\\&\leq&
(k+m+1)T(r,g)+S(r,g),
\end{eqnarray}
which is a contradiction.

{\bf Case 2}. $b=0$, $a\not=1$. From (\ref{leilee}), we have
\begin{equation}\label{leidedgagle}
F=\frac{G+a-1}{a}.
\end{equation}

Similarly, we also can get a contradiction, Thus, $a=1$ follows, it
implies that $F=G$.

{\bf Case 3}. $b=-1$, $a\not=-1$. From (\ref{leilee}), we obtain
\begin{equation}\label{leile}
F=\frac{a}{a+1-G}.
\end{equation}
Similarly, we can get a contradiction, $a=-1$ follows. Thus, we get
$F\cdot G=1$. From Lemma \ref{nan} and Lemma \ref{jintian}, we get
$f=tg$ for $t^{m}=t^{n+1}=1$. Thus, we get the proof of Theorem
\ref{2010d}.

\section{Discussions}

In this paper, we investigated the uniqueness of derivative of
difference polynomial of entire functions. It is an open question
under what conditions Theorem \ref{2010coro} holds for meromorphic
functions with $\rho_{2}(f)<1$. In addition, if
$[f^{n}(f^{m}-1)\Delta_{c}f]^{(k)}$ and
$[g^{n}(g^{m}-1)\Delta_{c}g]^{(k)}$ share one common value, we
believe that $f=tg$ for $t^{m}=t^{n+1}=1$. Unfortunately, we have
not succeed in proving that.

\label{lastpage-01}

\begin{thebibliography}{99}

\bibitem{ww} W. Bergweiler and A. Eremenko, \emph{On the singularities of the inverse to a meromorphic function
of finite order, Revista Matem\'{a}tica Iberoamericana.} \textbf{11}
(1995), 355--373.

\bibitem{chen} H. H. Chen and M. L. Fang , \emph{On the value distribution
of $f^{n}f'$}, Sci. China Ser. A. \textbf{38} (1995), 789--798.

\bibitem{chen22} H. H. Chen, \emph{Yoshida functions and Picard values of integral
functions and their derivatives,} Bull. Austral. Math. Soc.
\textbf{54} (1996), 373--381.

\bibitem{chiang} Y. M. Chiang and S. J. Feng, \emph{On the Nevanlinna
characteristic $f(z+\eta)$ and difference equations in the complex
plane,} The Ramanujan. J. \textbf{16} (2008), 105--129.

\bibitem{R.G} R. G. Halburd and R. J. Korhonen, \emph{Difference analogue of the lemma on the logarithmic derivative with
application to difference equations,} J. Math. Anal. Appl.
\textbf{314} (2006), 477--487.

\bibitem{tohge} R. G. Halburd, R. J. Korhonen and K. Tohge,
\emph{Holomorphic cures with shift-invariant hyperplane preimages},
arXiv: 0903-3236.

\bibitem{hal} R. G. Halburd and R. J. Korhonen, \emph{Meromorphic solutions of difference equations,
integrability and the discrete $Painlev\acute{e}$ equations}, J.
Phys. A. \textbf{40} (2007), 1--38.

\bibitem{hayman1} W. K. Hayman, \emph{Picard values of meromorphic functions
and their derivatives,} Ann. Math. \textbf{70} (1959), 9-42.

\bibitem{hayman3} W. K. Hayman, \textit{Meromorphic functions.} Oxford at the
Clarendon Press, 1964.

\bibitem{I.L} I. Laine, \textit{Nevanlinna Theory and Complex Differential
Equation.} Studies in Mathematics 15, Walter de Gruyter, Berlin-New
(1993).

\bibitem{I.L and yang} I. Laine and C. C. Yang, \emph{Value distribution of difference
polynomials,} Proc. Japan Acad. Ser. A \textbf{83} (2007), 148--151.

\bibitem{liukai arch} K. Liu and L. Z. Yang, \emph{Value distribution of the difference
operator}, Arch. Math. \textbf{92} (2009), 270--278.

\bibitem{rocky} K. Liu, \emph{Value distribution of differences of meromorphic
functions,} to appear in Rocky Mountain J. Math.

\bibitem{liucao} K. Liu, X. L Liu, and T. B Cao, \emph{Value distributions and
uniqueness of difference polynomials,} Advances in Difference
Equations Volume (2011), Article ID 234215, pp.12.

\bibitem{liuliucao} K. Liu, X. L Liu, and T. B Cao, \emph{Some results on zeros and uniqueness of difference-differential
polynomials}. Submitted.

\bibitem{liuzhangyang} K. Liu, C. H. Zhang, and L. Z. Yang, \emph{Uniqueness of entire functions and difference
polynomials}. Submitted.

\bibitem{lisheng} S. Li and Z. S. Gao, \emph{Finite order meromorphic
solutions of linear difference equations}, Proc. Japan Acad. Ser. A
\textbf{87} (2011), 73--76.

\bibitem{luolin}
X. D. Luo and W. C. Lin, \emph{Value sharing results for shifts of
meromorphic functions}, J. Math. Anal. Appl, \textbf{377} (2011)
441-449

\bibitem{mues} E. Mues, \emph{\"{U}ber ein Problem von Hayman,} Math.
Z, \textbf{164} (1979), 239--259.

\bibitem{Mohon} A. Z. Mohon'ho, \emph{The Nevanlinna characteristics of certain meromorphic functions}, Teor. Funktsii Funktsional. Anal. i Prilozhen.
\textbf{14} (1971), 83--87 (Russian).

\bibitem{Qi} X. G. Qi, L. Z. Yang and K. Liu, \emph{Uniqueness  and periodicity of meromorphic functions concerning
difference operator}, Comput. Math. Appl \textbf{60} (2010),
1739--1746.

\bibitem{wangyuefei} Y. F. Wang, \emph{On Mues conjecture and Picard values,} Sci. China.
\textbf{36} (1993), 28--35.

\bibitem{wangyuefei and Fang} Y. F. Wang and M. L. Fang, \emph{Picard values and normal families of
meromorphic functions with multiple zeros,} Acta Math. Sinica,
\textbf{14} (1) (1998), 17--26.

\bibitem{xu} J. F. Xu and H. X. Yi, \emph{Uniqueness of entire functions and differential polynomials}, Bull.
Korean Math. Soc. \textbf{44} (2007), 623--629.

\bibitem{hua} C. C. Yang and X. H. Hua, \emph{Uniqueness and value sharing of meromorphic
functions,} Ann. Acad. Sci. Fenn. Math. \textbf{22} (1997),
395--406.

\bibitem{yang and yi} C. C. Yang and H. X. Yi, \textit{Uniqueness Theory of Meromorphic
Functions.} Kluwer Academic Publishers (2003).

\bibitem{zhang11} J. L. Zhang \emph{Value distribution and shared sets of differences of meromorphic
functions}, J. Math. Anal. Appl, \textbf{367} (2010), 401--408.

\end{thebibliography}
\end{document}